**Note on Integer Factoring Methods IV**
**N. A. Carella, September, 2008.**

**Abstract:** This note continues the theoretical development of deterministic integer factorization algorithms based on systems of polynomials equations. This work exploits a new idea in the construction of irreducible polynomials with parametized roots, and recent advances in polynomial lattices reduction methods. The main result establishes a new deterministic time complexity bench mark in the theory of integer factorization.

## 1 Introduction

This note presents a deterministic integer factorization algorithm based on a system of polynomials equations. This technique combines a new irreducible polynomials construction technique and recent advances in lattice reduction methods to obtain a new result. The main result establishes a new deterministic time complexity bench mark. Background materials in the theory of integer factorization are given in [CE], [CP], [LA], [MZ], [RL], [S], [W], and similar sources.

The second section recalls the known results on the time complexity of integer factorization. It continues with the main contributions, Lemma 3 and Theorem 5, and concludes with an algorithm.

## 2 Main Contributions

This work builds on the earlier successful applications of the theory of polynomials equations and lattice reduction methods to integer factorization.

*Previous Results*
The previous works claim the followings.

***Theorem 1.*** ([CR]) If the $\log_2(N)/4$ least significant bits of a factor $p$ of $N$ are known, then the factorization of the integer $N = pq$, $p < q < 2p$, has deterministic logarithmic time complexity $O((\log N)^c)$, $c > 0$ constant.

In the Summer of 2007 this result was improved to the following.

***Theorem 2.*** Let $N = pq$, $p < q < 2p$. If the $(1/6)\log_2(N)$ most significant bits of a factor $p$ are known, then the factorization of $N$ has deterministic logarithmic time complexity $O((\log N)^c)$, $c > 0$ constant.

**Note.** The standard term *polynomial time* has been replaced with the more descriptive term *logarithmic time*. This is patterned after the closely related term *exponential time*.

*Construction of Irreducible Polynomials*

The height of a polynomial $f(x, y, z) = \sum_{0 \le i,j,l \le d} a_{i,j,l} x^i y^j z^l \in \mathbb{Z}[x,y,z]$ of maximum degree $\deg(f) = d$ in the variables $x$, $y$ and $z$ is given by the expression $\| f(x,y,z) \|_\infty = \max\{ \, |a_{i,j,k}| : 0 \le i,j,\, k \le d \, \}$.

**Lemma 3.** Let $\alpha, \beta > 0$ be a pair of parameters, and let $N = pq$ be a composite integer such that $p = O(N^\alpha)$ and $q = O(N^{1-\alpha})$. Then there exists an irreducible polynomial $f(x, y, z) = c_4 xy + c_3 x + c_2 y + c_1 z + c_0 \in \mathbb{Z}[x, y, z]$ with the following properties.

(i) The polynomial $f(x, y, z \,)$ has a small integer root $(x_0, y_0, z_0)$ where $0 \le |x_0| \le X \le N^\alpha$, $0 \le |y_0| \le Y \le N^{1-\alpha}$, and $0 \le |z_0| \le Z \le O((\log N)^B)$, $B > 0$ constant.
(ii) The factors of $N$ can be written as $p = mx_0 + c$ and $q = ny_0 + d$, where the moduli $n$ and $m$ (possibly relatively prime) are of size $O((\log N)^A)$, $A > 0$ constant.
(iii) The height of $f(xX, yY, zZ)$ satisfies the inequality $\| f(xX, yY, zZ) \|_\infty \ge N^{1+\beta}$.
(iv) The polynomial can be generated in deterministic logarithmic time $O((\log N)^c)$, $c > 0$ constant.

Proof: Let $n$ and $m$ be (possibly relatively prime) moduli of sizes $O((\log N)^A)$, and let $k = O(N^\beta)$ be an integer, with $A > 0$ and $\beta > 0$ constants. Next rewrite the equation $f(x, y) = xy - N$ as an equation of three variables

$$f(x, y, z) = (mx + c)(ny + d)(kz + e) - rN = c_4 xy + c_3 x + c_2 y + c_1 z + c_0, \qquad (1)$$

where $1 \le c, d < O((\log N)^A)$, and $r = kz_0 + e$ is prime (or nearly prime) with $|z_0| \le Z \le O((\log N)^B)$. The coefficients $c_i$ are obtained after a selective replacement of the known variable $z = z_0$. Clearly this is an irreducible polynomial over the integers and of height $\| f(xX, yY, zZ) \|_\infty \ge rN = N^{1+\beta}$, and has a small integral root $0 \le |x_0| < N^\alpha$, $0 \le |y_0| < N^{1-\alpha}$, $0 \le |z_0| < O((\log N)^B)$. ∎

Apparently, it is a very difficult proof for those that have not seem it before. Nevertheless, it is an elementary transformation/deformation of the most important polynomial $f(x, y) = xy - N$ in integer factorization.

The basic algorithm of Lemma 3 is sketched below, it is designed to work in tandem with Algorithm II, also note that the data $m, n, c, d$ can be either inputted or internally generated.

*Algorithm* I
Input: $\alpha, \beta > 0$, and $N = pq$ such that $p = O(N^\alpha)$.
Output: $f(x, y, z) = c_4 xy + c_3 x + c_2 y + c_1 z + c_0$ and $m, n, c, d$.
1. Set $T = (\log N)^A$, $A > 1$, and generate a pair of primes or nearly primes moduli $n$ and $m < T$.
2. Generate a prime $r = kz_0 + e$ (or nearly prime) with $k = O(N^\beta)$ and $|z_0| \le Z \le O((\log N)^B)$.
3. Select a pair $c < m$, and $d < n$.
4. Compute the coefficients of the irreducible polynomial $f_{c,d}(x, y, z) = c_4 xy + c_3 x + c_2 y + c_1 z + c_0$.
5. Return $f(x, y, z) = f_{c,d}(x, y, z)$.

*Algebraically Independent Polynomials*
Although the technique of Lemma 3 can generate one or more irreducible polynomials, these polynomials are not algebraically independent. Accordingly, lattice reduction method is utilized to generate another algebraically independent polynomial. Further, since the third variable is known or its value is very small and can be determined by brute force search, just one additional algebraically independent polynomial is required.

**Theorem 4.** ([ER]) Let $f(x, y, z) = c_4xy + c_3x + c_2y + c_1z + c_0 \in \mathbb{Z}[x, y, z]$ be an irreducible polynomial of height $\| f(xX, yY, zZ) \|_\infty = W$ and with a small integer root $(x_0, y_0, z_0)$ such that $| x_0 | < X$, $| y_0 | < Y$ and $| z_0 | < Z$. Suppose that the inequality

$$X^{3+3\tau} Y^{3+6\tau+3\tau^2} Z^{2+3\tau} < W^{2+3\tau-\varepsilon}, \tag{2}$$

where $\tau > 0$ is a lattice parameter, holds. Then there exists a pair of linearly independent polynomials $f_1(x, y, z)$ and $f_2(x, y, z)$ not multiple of $f(x, y, z)$, with a common root. Furthermore, the polynomials are generated in deterministic logarithm time.

The complete analysis of this and other special cases of polynomials in three and four variables and the corresponding polynomials bases of the polynomials lattices are given in [ER], and [JM].

The two polynomials $f_1(x, y, z)$ and $f_2(x, y, z)$ are linearly independent, but not guaranteed to be algebraically independent. However, the two pairs of polynomials $f(x, y, z), f_1(x, y, z)$ and $f(x, y, z), f_2(x, y, z)$ are guaranteed to be algebraically independent. Recent advances in the construction of three algebraically independent polynomials are discussed in [BA].

*The Main Result*
In the last decades the techniques of the theory of polynomials equations and lattice reduction methods have emerged as powerful tools in the theory of integer factorization.

**Theorem 5.** The factorization of a composite integer $N \in \mathbb{N}$ has deterministic logarithmic time complexity $O((\log N)^c)$, $c > 0$ constant.

Proof: Without loss in generality, let $N = pq$ be a balanced integer, $p < q < 2p$. Put $\alpha = 1/2$, and let $\beta = 1/2 + \gamma$ for some $\gamma > 0$, and fix a pair of moduli $n, m = O((\log N)^A)$, where $A > 0$ is a constant. Then it is clear that the integer $N$ has its factors in some residue classes

$$p = mx + c \quad \text{and} \quad q = ny + d, \tag{3}$$

where $0 \leq | c |, | d | < O((\log N)^A)$. At most $O((\log N)^{2A})$ pairs $(c, d)$ has to be tested to determine the correct residues classes (3) of the factors. Given the correct pair $(c, d)$, there exists an irreducible polynomial

$$f(x, y, z) = c_4xy + c_3x + c_2y + c_1z + c_0 \tag{4}$$

over the integers $\mathbb{Z}$, which has a small integer solution $(x_0, y_0, z_0)$ such that

$$p = mx_0 + c \quad \text{and} \quad q = ny_0 + d, \tag{5}$$

where $0 \leq | x_0 |, | y_0 | < N^{1/2} = N^\alpha$ and $0 \leq | z_0 | \leq O((\log N)^B)$, $B > 0$ constant, see Lemma 3. Moreover, the height satisfies the inequality $W = \| f(xX, yY, zZ)) \|_\infty \geq N^{1+\beta} \geq N^{3/2+\gamma}$. Now by Theorem 4, there exists another algebraically independent polynomial $g(x, y, z)$ that shares the same root $(x_0, y_0, z_0)$ and it can be determined using lattice reduction techniques whenever the inequality

$$X^{3+3\tau} Y^{3+6\tau+3\tau^2} Z^{2+3\tau} < W^{2+3\tau-\varepsilon}, \tag{6}$$

holds. Replacing $X < N^{1/2}$, $Y < N^{1/2}$, $Z < O(N^\delta)$ and $N^{3/2+\gamma} \leq W$ in (6) returns

$$X^{3+3\tau}Y^{3+6\tau+3\tau^2}Z^{2+3\tau} < N^{3+9\tau/2+3\tau^2/2+(2+3\tau)\delta} \leq N^{(3/2+\gamma)(2+3\tau-\varepsilon)} \leq W^{2+3\tau-\varepsilon}, \tag{7}$$

where $\delta > 0$ is an arbitrarily small number, and all the relevant constants has been omitted. These data in turn imply that (6) and (7) holds if and only if

$$\frac{3\tau^2/2+(2+3\tau)\delta+3\varepsilon/2}{2+3\tau-\varepsilon} < \gamma \tag{8}$$

holds. Further, since there is almost no restriction on the parameter $\gamma > 0$, for example, $\gamma = 1/4$ or $1/3$ or $1/2$, etc is feasible, the previous inequalities (6), (7) and (8) hold for any appropriate choice of $\gamma$.

Ergo the solution $(x, y, z) = (x_0, y_0, z_0)$ of the system of equations

$$f(x, y, z) = 0, \; g(x, y, z) = 0,$$

can be recovered by means of resultants or Grobner bases calculations. Specifically, computing the roots of the polynomials

$$r_1(x) = \text{Res}_y(f(x, y, z_0), g(x, y, z_0)) \quad \text{and} \quad r_2(y) = \text{Res}_x(f(x, y, z_0), g(x, y, z_0)), \tag{9}$$

where $z_0$ is known. Next observe that and each of these algorithms above has logarithmic time complexity. In particular, the running time of the entire algorithm is dominated by at most $O((\log N)^{2A})$ lattice reduction steps, one for each pairs $(c, d)$. Thus, the overall time complexity of the integer factorization algorithm is deterministic logarithmic time $O((\log N)^c)$, $c > 0$  constant.                    Quod erat demonstrandum    ∎

The choice of parameter $\alpha > 0$ assumes a priori knowledge on the sizes of the factors $p = N^\alpha$ and $q = N^{1-\alpha}$ of $N = pq$. Furthermore, since the subset of balanced integers $\mathcal{B} = \{ N = pq : p < q < ap$, with $p, q$ primes $\}$ is the most important case in integer factorization, it was set to $\alpha = 1/2$. Balanced integers are the most difficult to factor. However, the probability of an arbitrary integer of being balanced is negligible. Indeed, the subset of balanced integers has zero density in the set of integers. More precisely, it has cardinality $\mathcal{B}(x) = \#\{ N = pq \leq x : p < q < ap \} = c_0 x/\log(x)^2$, $c_0 = c_0(a)$ constant, see [DM].

Numerical experiments will have to be performed to determine the best choices of the parameters $\beta = 1/2 + \gamma$ and $\tau > 0$. The first controls the height of the polynomial $f(x, y, z)$ and the second is part of the lattice basis, see Theorem 4.

The last algorithm below encodes the basic procedure of Theorem 5. In step 2.1, the data $m, n, c, d$ is passed on to Algorithm I.

*Algorithm* II
Input: $N = pq$, and $\alpha, \beta > 0$ such that $p = O(N^\alpha)$.
Output: $p, q$.
1. Set $T = (\log N)^A$, $A > 0$, and select a pair of primes or nearly prime moduli $n$ and $m < T$.
2. For $c, d < T$ do
2.1 Construct an irreducible polynomial $f_{c,d}(x, y, z) = c_4 xy + c_3 x + c_2 y + c_1 z + c_0$ such that  $0 \leq |x_0| \leq N^\alpha$, $0 \leq |y_0| \leq N^{1-\alpha}$, and $0 \leq |z_0| \leq O((\log N)^B)$, $B > 0$ constant, see Lemma 3 and Algorithm I.
2.2 Construct an algebraically independent polynomial $g_{c,d}(x, y, z)$, see Theorem 4.
2.3 Compute the root $(x_0, y_0, z_0)$ of the system of equations $f(x, y, z) = 0$, $g(x, y, z) = 0$, using resultants or Groebner bases methods, here $z_0$ is known.

2.4 Compute the potential factors $p_{c,d} = mx_0 + c$ and $q_{c,d} = ny_0 + d$.

2.5. If $1 < \gcd(p_{c,d}, N) < N$ or $1 < \gcd(q_{c,d}, N) < N$, then halt.

3. Return $p = p_{c,d}$, and $q = N/p_{c,d}$ or $q = q_{c,d}$, and $p = N/q_{c,d}$.

Ultimately an algorithm that accepts an arbitrary integer and internally generates all its parameters seems to be feasible in the near future.

## 3 Polynomials Equations

Although the problem on hand is effective integer factorization, the solution proposed here is almost entirely based on the analytic and algebraic theory of polynomials, polynomials equations, and related topics. The number theoretical aspect of this problem is almost completely absent in the proposed solution.

A short introduction to polynomials and systems of polynomials equations is supplied in this section. The reader should consult the literature to bridge the gaps, some standard references are [CO], [MT], [PV], [RD], [MR], [RS], [SL], and so on.

### 3.1 Univariate Polynomials

A polynomial equation $f(x) = 0$ is the simplest system of polynomials equation, and successive applications of resultants or Groebner bases techniques can reduce a system of polynomials equations

$$f_1(x_1, \ldots, x_n) = 0, \ldots, f_n(x_1, \ldots, x_n) = 0$$

to an ordered system of polynomials equations

$$g_1(x_1, \ldots, x_n) = 0, \ldots, g_{n-1}(x_1, x_2) = 0, g_n(x_1) = 0,$$

which involves the determination of the roots of a univariate polynomial. In light of this, it is natural that certain aspects of the theory of univariate polynomials are crucial in the development of root finding algorithms to solve systems of polynomials equations. The numerical analysis of roots finding algorithms are not covered here, see [PT] and similar sources.

**Theorem 6.** (Argand 1806) A polynomial $f(x) \in \mathbb{C}[x]$ of degree $n = \deg(f)$ has $n$ complex zeroes in the field of complex numbers $\mathbb{C}$.

Previous works on the real roots were stated by Girard in 1629, and other authors, and are many proofs of this result, including two by Gauss are known. Essentially the same claim holds with the field of complex numbers replaced by other fields or rings mutatis mutandis.

**Theorem 7.** The roots of a polynomial are expressible in terms of elementary operations and radicals if and only if the Galois group of the polynomial is solvable.

A group is *solvable* if the sequence $\{ e \} = G_v \subset G_{v-1} \subset \cdots \subset G_1 \subset G_0 = G$ has abelian quotient $G_{i-1} / G_i$ and $G_i$ is normal in $G_{i-1}$.

This result combines the works of Ruffini 1799, Abel 1824, Galois 1832, and other contributors. It applies to every polynomial of degree less than five and many other polynomials of higher degrees. It is quite simple to compute a polynomials that does not have a solvable Galois group, exempli gratia, the roots of $x^5 - x - 1$, $x^5 + x - 1$, $x^5 + x + 1$, $x^5 - 4x + 2$, $\ldots$, are not given in terms of radicals and elementary arithmetic operations since the Galois group of each of this polynomial is the symmetric group $S_5$.

**Height Of Polynomials**

Let $|x|_v$ be an absolute value. The height of the polynomial $f(x) = a_nx^n + \cdots + a_1x + a_0$ is defined by, and let $H(f)$ = max $\{\,|a_i|_v\,\}$. The measure $H(f)$ is an extension of the height $\|\,r/s\,\|_\infty = \max\{\,|\,r\,|,\,|\,s\,|\,\}$ of a rational number $r/s$, both notations $H(f) = \|\,f\,\|_\infty = \max\,\{\,|a_i|_v\,\}$ are widely used. The Lehmer measure of $f(x)$ is defined by

$M(f) = |a_n| \prod_{i=1}^{n} \max\{1, |\alpha_i|\}$. This metric was first used by Lehmer to investigate the prime values of

polynomials. The Lehmer's problem is concerned with the determination of the smallest measure $M(f) > 1$, the polynomial $x^{10} + x^9 - x^7 - x^6 - x^5 - x^4 - x^3 + x + 1$ of measure $M(f) = 1.1762808\ldots$ holds the record.

The height accounts for the extreme coefficients, and the Lehmer measure accounts for the large roots $|\,z\,| > 1$ outside the unit disk. For example, a product of cyclotomic polynomials has a Lehmer measure of $M(f) = 1$, otherwise it seems to be $M(f) > 1$.

**Ranges of the Roots**

***Theorem* 8.** (Cauchy 1829)  Let $\theta$ be a root of $f(x) = a_nx^n + \cdots + a_1x + a_0 \in \mathbb{C}[x]$. Then $|\,\theta\,| < (1 + \|\,f\,\|_\infty)/|\,a_n\,|$.

Proof: Rearrange the equation $f(x) = 0$ and take absolute values and height:

$$\left|\,a_n\theta^n\,\right| = \left|\,a_{n-1}\theta^{n-1} + \cdots + a_0\,\right| \leq \left|\,\theta^{n-1} + \cdots + \theta + 1\,\right| \|f\|_\infty \leq \frac{|\,\theta\,|^n - 1}{|\,\theta\,| - 1} \|f\|_\infty. \qquad (10)$$

Again rearrange the both side to arrive at the claim. ∎

In short, a polynomial has all its roots $z_i$ in within a disk of radius $(1 + \|\,f\,\|_\infty)/a_n$. The auxiliary polynomial $g(x) = |\,a_n\,|x^n + \cdots + |\,a_1\,|x + |\,a_0\,| \in \mathbb{R}[x]$ has a maximal real root $|\,r_0\,|$ such that $|\,z_i\,| \leq |\,r_0\,|$, $1 \leq i \leq n.$, and it can be computed using numerical algorithms.

***Theorem* 9.** (Enestrom 1893)  Let $\theta$ be a root of $f(x) = a_nx^n + \cdots + a_1x + a_0 \in \mathbb{C}[x]$ with $a_n \geq 0$. Then

$$\min\left\{\frac{a_i}{a_{i-1}} : 1 \leq i < n\right\} \leq \left|\,\theta\,\right| \leq \max\left\{\frac{a_i}{a_{i-1}} : 1 \leq i \leq n\right\}. \qquad (11)$$

The reader should confer the literature for refined and recent versions of these estimates.

**Polynomials Expansions**

The expansion $f(x) = a_nx^n + \cdots + a_1x + a_0 = b_nL_n(x) + \cdots + b_1L_1(x) + b_0L_0(x)$, $b_i \in \mathbb{C}$, of a polynomial in terms of a basis $\{\,L_n(x),\,\ldots,\,L_0(x)\,\} \subset \mathbb{C}[x]$ can often reveal special properties of the roots of $f(x)$ and other spectral information. A demonstration of this concept is the next result.

***Theorem* 10.** (Turan 1950)  Let $f(x) = b_nH_n(x) + \cdots + b_1H_1(x) + b_0H_0(x) \in \mathbb{R}[x]$ be the Hermite expansion of $f(x)$. If the inequality $\sum_{0 \leq k \leq n-2} 2^k k! b_k^2 < 2^2(n-1)! b_n^2$ holds, then $f(x)$ has $n$ distinct real roots.

The factorial and Legendre expansions of a polynomial are respectively given by

$$f(x) = \sum_{k=0}^{n} b_k \binom{x}{k} \quad \text{and} \quad f(x) = \sum_{i=0}^{n} y_i L_i(x), \text{ where } F_i(x) = \binom{x}{i} \quad \text{and} \quad L_i(x) = \prod_{i \neq j}\left(\frac{x - x_j}{x_i - x_j}\right) \qquad (12)$$

are the bases and $(x_0, y_0)$, …, $(x_n, y_n)$ are the preassigned values $f(x_i) = y_i$ of $f(x)$.

**Number Of Real Roots**

***Theorem 11.*** (Kac 1943?)  If the coefficients of the polynomial $f(x) = a_n x^n + \cdots + a_1 x + a_0 \in \mathbb{R}[x]$ are normally distributed with mean $\mu = 0$ and variance $\sigma^2 > 0$, then the expected number $E_n(f)$ of real zeros of $f(x)$ is asymptotic to $(2/\pi)\log n$.

There is an exact formula for $E_n(f) = (2/\pi)\log n + o(\log n)$. Moreover, about the same asymptotic number of real roots of $f(x)$ also hold for other distributions of the coefficients of $f(x)$, see the literature.

***Theorem 12.*** ( Erdos-Turan 1950)  Let $n_f[\alpha, \beta]$ denotes the number of zeroes $z = re^{i\theta}$ of $f(x) = a_n x^n + \cdots + a_1 x + a_0$ with $\theta \in [\alpha, \beta] \subset [0, 2\pi]$. Then

$$\left| \frac{n_f[\alpha, \beta]}{n} - \frac{\beta - \alpha}{n} \right|^2 \leq 2 \log\left( \frac{|a_0| + |a_1| + \cdots + |a_n|}{\sqrt{|a_0 a_n|}} \right). \tag{13}$$

This result is a multitask tool: it sheds lights on the distribution of the zeroes on a suitable disk and provides an estimate of the maximum number of positive roots (or negative) of a polynomial as a function of its coefficients, et cetera.

It is quite easy to derive the inequality

$$R^2 \leq 2n \log\left( \frac{|a_0| + |a_1| + \cdots + |a_n|}{\sqrt{|a_0 a_n|}} \right), \tag{14}$$

for the number $R$ of positive (or negative) zeroes of $f(x) = a_n x^n + \cdots + a_1 x + a_0 \in \mathbb{R}[x]$ with $n^{-t} \leq |a_n| \leq n^t$, some $t > 0$.

The constraint $n^{-t} \leq |a_n| \leq n^t$, some $t > 0$, on the coefficients is crucial. For the weaker $0 \leq |a_n| \leq n^t$ constraint, the polynomials $f(x) = (x^{2k} - a_1)(x^{2k} - a_2) \cdots (x^{2k} - a_m)$ with small $a_i > 0$, in $\mathbb{R}[x^2]$ are counterexamples. Each of these polynomials has $2m$ real roots almost completely independent of the degree $2km = \deg(f(x))$ of $f(x)$.

A *Sturm sequence* $\{ f_0(x), f_1(x), …, f_{k+1}(x) \}$ is defined by the intermediate results of the Euclidean algorithm $f_{i-1}(x) = q_i(x) f_i(x) - f_{i+1}(x)$, with the initial conditions $f_0(x) = f(x), f_1(x) = f'(x)$.

***Theorem 13.*** (Sturm 1835)  Let $\{ f_0(x), f_1(x), …, f_{k+1}(x) \}$ be the Sturm sequence of $f(x)$ in $\mathbb{R}[x]$. Then

$$Z(f_0(\alpha) f_1(\alpha), …, f_{k+1}(\alpha)) - Z(f_0(\beta) f_1(\beta), …, f_{k+1}(\beta)),$$

where $Z(x_0, x_1, …, x_m)$ is the number of signs changes in the sequence $\{ x_0, x_1, …, x_m \}$, is the number of real roots of $f(x)$ in the interval $[\alpha, \beta]$.

The proofs of several versions and generalizations of this result are explored in [PV], and the algorithm is analyzed in [BW, p. 400].

The center of mass of the roots of a polynomial is the expression $\hat{z} = z - n\dfrac{f(z)}{f'(z)}$. It is a sort of a geometric average of the roots of $f(x)$ with respect to a complex number $z \in \mathbb{C}$.

**Theorem** 14. (Laguerre 1860?) Let $z \in \mathbb{C} - \mathbb{R}$. Then the roots of $f(z)$ are all real numbers if and only if $(\operatorname{Im} z)(\operatorname{Im} \hat{z}) < 0$.

The proof appears in the literature, [Am. Math, Month, Vol. 86, No. 8, p. 648-658, 1979]. The choice $z = i$ realizes a very simple *all real roots test*, which is far more effective than the discriminant based test for polynomial of degree $> 2$.

Let $z_1, z_2, \ldots, z_n$, be the roots of a polynomial $f(x)$ of degree $\deg(f(x)) = n$. The separation or distance between a pair of roots is the real number $|z_i - z_j|$.

The expected distance $E(|z - z_i|)$ of the roots $z_i$ of a polynomial from a complex number $z \in \mathbb{C}$ is $|f(z)|^{1/n} > 0$. This is obtained by means of the arithmetic-geometric mean inequality

$$\frac{x_1 + x_2 + \cdots + x_n}{n} \geq (x_1 \cdot x_2 \cdots x_n)^n. \tag{15}$$

A more precise result is the average separation of the roots $|z_i - z_j|$, which is also proved using the geometric-arithmetic mean inequality.

**Theorem** 15. The average separation of the roots $z_i$ of a squared free polynomial is $E(|z_i - z_j|) \geq \left| \operatorname{Re} s(f, f', x) \right|^{1/n(n-1)} > 0$.

A root near the complex number $z \in \mathbb{C}$ is computed using a variety of algorithms, [PT]. For example, the Newton's method for squared free polynomial generates a sequence of numbers $z_{n+1} = z_n - \dfrac{f(z_n)}{f'(z_n)}$, $n \geq 0$, converging to a root $\theta$ of $f(z)$. The sequence is uniformly convergent, and the speed of convergence is quadratic.

The rational roots $z = r/s$ of a polynomial $f(x) = a_n x^n + \cdots + a_1 x + a_0 \in \mathbb{Z}[x]$ are often found by trying the factors of the last and first coefficients, that is, $r \mid a_0$ and $s \mid a_n$. But there are two obstacles whenever these coefficients are large integers:
(i) Obtaining the factorizations of $a_0$ and $a_n$. (ii) The possibility of exponentially many factors in the factorizations of $a_0$ and $a_n$.

## Magnitudes Of The Coefficients

Write a polynomial in the form $f(x) = a_n x^n + \cdots + a_1 x + a_0 = a_n \prod_{i=1}^{n} (x - z_i)$. Expanding and rewriting it gives

$f(x) = a_n(x^n - \sigma_1 x^{n-1} + \sigma_2 x^{n-2} + \cdots + \sigma_{n-1} x + (-1)^n \sigma_n)$, where

$$\sigma_1 = z_1 + \cdots + z_n, \ldots, \sigma_k = \sum_{1 \leq i_1 < i_2 < \cdots < i_k \leq n} z_{i_1} z_{i_2} \cdots z_{i_k}, \ldots, \sigma_n = z_1 \cdots z_n, \tag{16}$$

are the symmetric functions. The power functions are defined by $\rho_k = z_1^k + \cdots + z_n^k$. The generating functions of the symmetric and power functions are respectively

$$F(T) = \prod_{i=1}^{n}(1 + z_i T) \quad \text{and} \quad G(T) = \sum_{i=1}^{n}\frac{z_i}{1 - z_i T}. \tag{17}$$

The corresponding Newton's recursive formula is $n\sigma_n = \sum_{k=1}^{n}(-1)^{k-1}\rho_k \sigma_{n-k}$, which is a sort of convolution of power and symmetric functions.

The estimate for the coefficients of a factor $g(x)$ of a polynomial $f(x)$ is essentially patterned after the extreme polynomial $f(x) = (x + c)^n$. Clearly, a coefficients of a factor $g(x) = (x + c)^m = b_m x^m + \cdots + b_1 x + b_0$ satisfies the expression

$$\mid b_k \mid = \binom{m}{k}c^k \leq \mid a_k \mid = \binom{n}{k}c^k. \tag{18}$$

**Proposition 16.** ([MT]) Let $z_1 \geq 1$, $z_2 \geq 1$, ..., $z_n \geq 1$ be real numbers, and let $\sigma_k = \sum_{1 \leq i_1 < i_2 < \cdots < i_k \leq n} z_{i_1} z_{i_2} \cdots z_{i_k}$ be the $k$th symmetric sum in the variables $z_1$, $z_2$, ..., $z_n$. Then

$$\mid \sigma_k \mid \leq \binom{n-1}{k-1}\mid z_1 z_2 \cdots z_n \mid + \binom{n-1}{k}. \tag{19}$$

**Proposition 17.** ([MT]) If $g(x) = b_m x^m + \cdots + b_1 x + b_0$ divides $f(x) = a_n x^n + \cdots + a_1 x + a_0$ then

$$\text{(i)} \ \mid b_k \mid \leq \binom{m}{k}\parallel f \parallel_\infty \ \text{and (ii)} \ \mid b_k \mid \leq \binom{n-1}{k}\parallel f \parallel_2 + \binom{n-1}{k-1}\mid a_m \mid. \tag{20}$$

**Theorem 18.** (i) $\mid a_k \mid \leq \binom{n}{k}M(f)$ (ii) $H(f) \leq 2^{n-1}M(f)$.

**Lemma 19.** (Gauss) If $\mid x \mid_v$ is a nonarchemidean absolute value, then $H(fg) = H(f)H(g)$.

This is a modern version of the Gauss's lemma, the original proof was based on the gcd's of the coefficients of the polynomials $f(x)$ and $g(x)$, the proof is given in [PV].

## Irreducibility and Factorization

A commutative ring $\mathbf{R}$ is a *unique factorization domain* if each element $x \in \mathbf{R} - \mathbf{R}^*$ has a unique decomposition as a product $r = r_1^{e_1} r_2^{e_2} \cdots r_v^{e_v}$, $r_i$ irreducible, up to a permutation and a unit. The units in the ring are the subset of invertible elements $\mathbf{R}^* = \{ x \in \mathbf{R} : x^{-1} \text{ exists} \}$.

**Example 20.** The rings of integers and polynomials $\mathbb{Z}$, $\mathbb{Z}[\sqrt{-1}]$, $\mathbb{Z}[\sqrt{2}]$, $\mathbb{Z}[\sqrt{3}]$, ..., $\mathbb{Z}[x]$, $\mathbb{Z}[\sqrt{-1}][x]$, $\mathbb{Z}[\sqrt{2}][x]$, $\mathbb{Z}[\sqrt{3}][x]$, ..., are UFDs, but rings of integers and polynomials $\mathbb{Z}[\sqrt{-5}]$, $\mathbb{Z}[\sqrt{-6}]$, $\mathbb{Z}[\sqrt{10}]$, $\mathbb{Z}[\sqrt{15}]$, ..., $\mathbb{Z}[\sqrt{-5}][x]$, $\mathbb{Z}[\sqrt{-6}][x]$, $\mathbb{Z}[\sqrt{10}][x]$, $\mathbb{Z}[\sqrt{15}][x]$, ..., are not UFDs.

The absence of a property such as integral closure in a ring disables the UFD property in some rings, (I do not know if this is true in any ring). For example, the quadratic ring $\mathbb{Z}[\sqrt{d}\,]$ is integrally closed if $d \not\equiv 1 \bmod 4$, otherwise, $\mathbb{Z}[(1+\sqrt{d})/2]$ is integrally closed if $d \equiv 1 \bmod 4$ but $\mathbb{Z}[\sqrt{d}\,]$ is not an UFD. Due to this difficulty, $\mathbb{Z}[\sqrt{13}\,]$ is not a UFD, in fact, $-20 = -4 \cdot 5 = (3-\sqrt{13})(3+\sqrt{13}) \cdot 5$, but $\mathbb{Z}[(1+\sqrt{13})/2]$ is an UDF, in this ring $-20 = -4 \cdot 5$.

**Theorem 21.** (Gauss 1799) If **R** is a unique factorization domain, then the ring of polynomials $\mathbf{R}[x_1, ..., x_n]$ is a unique factorization domain, $n \geq 1$.

A polynomial $f(x)$ is *irreducible* if $f(x) = g(x)h(x)$ implies that either $g(x) = $ constant or $h(x) = $ constant. Otherwise $f(x)$ is reducible. In some polynomials rings each polynomial has a unique factorization as a product of irreducible polynomials. There are many polynomials irreducibility tests, one of them is recorded here as an illustration and a reference.

**Theorem 22.** ([RM]) Let $f(x) = a_n x^n + \cdots + a_1 x + a_0$ be a polynomial of degree $n$ in $\mathbb{Z}[x]$ and set $H(f) = \max\{a_i/a_n : 0 \leq i \leq n-1\}$. If $f(n)$ is prime for some $n \geq H(f) + 2$ then $f(x)$ is irreducible in $\mathbb{Z}[x]$.

This is a general purpose test, it tests any polynomial independently of its form, but polynomials with small coefficients require more careful analysis. Observe that a reducible polynomial $f(x)$ in $\mathbb{Z}[x]$ can assume prime values, but at most $\deg f$ times. It is possible that a reducible polynomial $f(x)$ in $\mathbf{R}[x]$ can assume prime values infinitely often in some ring **R**?

**Theorem 23.** ([L]) An univariate polynomial with coefficients in an infinite UFD can be factored in logarithmic time .

The time complexity of polynomial factorization is measured as a power of $\log \| f(x) \|$. Indirectly, this is a claim on the time complexity of determining all the integral roots of a polynomial $f(x) \in \mathbf{R}[x]$, where **R** is an infinite ring. However, the authors claim that the algorithm is not practical.

**Theorem 24.** If **R** supports an Euclidean algorithm, then the ring of polynomials $\mathbf{R}[x_1, ..., x_n]$ supports an Euclidean algorithm.

A discussion of these ideas and a proof appears in [BW, p. 98].

### 3.2 Multivariate Polynomials

Let **R** be a commutative ring and let $\mathbf{R}[x_1, ..., x_n] = \{ f(x_1, ..., x_n) = \sum_{\alpha = (\alpha_1, ..., \alpha_n)} a_\alpha x_1^{\alpha_1} x_2^{\alpha_2} \cdots x_n^{\alpha_n} : \alpha_i \geq 0, a_\alpha \in \mathbf{R} \}$ be the ring of polynomial functions of $n$ variables.

The *index* set of the polynomial $f(x) = f(x_1, ..., x_n)$ is the subset $\mathrm{ind}(f) = \{ (\alpha_1, ..., \alpha_n) : 0 \neq \alpha_i \in \mathbb{N} \}$ of positive integers vectors $\alpha = (\alpha_1, ..., \alpha_n)$. A polynomial is often written as $f(x) = \sum_\alpha a_\alpha x^\alpha$ .

The *maximal* and the *total* degree of a polynomial of one or more variables are defined by the formulae $\deg_{\max}(f) = \max\{ \alpha_i : 1 \leq i < n, \alpha \in \mathrm{ind}(f) \}$ and $\deg_{\mathrm{total}}(f) = \max\{ \alpha_1 + \cdots + \alpha_n : \alpha \in \mathrm{ind}(f) \}$ respectively.

***Proposition* 25.** The linear space $L_d(x_1, ..., x_n) = \{ \ f(x) = \sum_{|\alpha| \le d} a_\alpha x^\alpha \in \mathbf{R}[x_1, ..., x_n] \ \}$ of monomials of total degree $|\alpha| = \alpha_1 + \cdots + \alpha_n \le d$ has dimension $\binom{n+d}{n} = \dfrac{(n+d)!}{d! \, n!}$.

## Division Algorithm

The division algorithm in a polynomials ring $\mathbf{R}[x_1, ..., x_n]$ states the followings: Given a fixed subset $\{ f_1(x), \ldots, f_m(x) \}$ and an ordering on the integers lattice $\mathbb{N}^n$, then every polynomial $f(x) \in \mathbf{R}[x_1, ..., x_n]$ has a representation as $f(x) = a_1 f_1(x) + \cdots + a_m f_m(x) + r(x)$, where $a_i(x), r(x) \in \mathbf{R}[x_1, ..., x_n]$.

***Definition* 26.** Let $\{ f_1, \ldots, f_m \} \subset \mathbf{R}[x_1, ..., x_n]$, an *ideal* in $\mathbf{R}[x_1, ..., x_n]$ is the subset

$$I(f_1, ..., f_m) = \{ \ f(x) = a_1 f_1 + \cdots + a_m f_m : a_i \in \mathbf{R}[x_1, ..., x_n] \ \}. \tag{21}$$

The ideal $I = I(f_1, \ldots, f_m)$ can be viewed as the *nonlinear* span of the polynomials basis $\{ f_1, \ldots, f_m \}$ over $\mathbf{R}[x_1, ..., x_n]$. This is often denoted as $I = <f_1, \ldots, f_m>$. The division algorithm implies that every polynomial $f(x) \in I(f_1, ..., f_m)$ has a representation as $f(x) = a_1 f_1(x) + \cdots + a_m f_m(x)$, where $a_i(x) \in \mathbf{R}[x_1, ..., x_n]$.

An ideal in $\mathbf{R}[x_1, ..., x_n]$ satisfies the followings properties.
(i) $0 \in I$,     (ii) $f, g \in I$ implies that $f + g \in I$,
(iii) $f \in I$ implies that $sf \in I$ for all $s \in \mathbf{R}[x_1, ..., x_n]$.

The *radical* $\sqrt{I}$ of an ideal $I$ in a ring is the subset $\sqrt{I} = \{ f \in \mathbf{R} : f^n \in I, \text{ some integer } n \ge 1 \}$. If $\sqrt{I} = I$, then $I$ is called a radical ideal. The definition is applicable to any ring $\mathbf{R}$. For example, $I = < 3 \cdot 5 > \subset \mathbb{Z}$ is a radical ideal, but $I = < 3^2 \cdot 5 > \subset \mathbb{Z}$ is not a radical ideal.

***Example* 27.** The ideal $I = < \sigma_1, \ldots, \sigma_n > \subset \mathbf{R}[x_1, ..., x_n]$ is a maximal ideal, in fact $I = \mathbf{R}[\sigma_1, ..., \sigma_n] = \mathbf{R}[\rho_1, ..., \rho_n]$, where $\rho_k$ and $\sigma_k$ are the power and symmetric functions of the variables $x_1, \ldots, x_n$. The verification follows from the observation that the system of equations $\sigma_1 = 0, \ldots, \sigma_n = 0$ has a unique solution $(x_1, \ldots, x_n) = (z_1, \ldots, z_n)$, so $I = < x_1 - z_1, \ldots, x_n - z_n >$.

***Theorem* 28.** (Hilbert 1888?)  Every ideal $I$ in $\mathbf{R}[x_1, ..., x_n]$ has a finite basis.

## 3.1 Systems Of Polynomial Equations

A system of polynomials equation is defined by $f_1(x_1, \ldots, x_n) = 0, \ldots, f_m(x_1, \ldots, x_n) = 0$.

The corresponding solution set

$$V(f_1, \ldots, f_m) = \{ \ (x_1, \ldots, x_n) \in \mathbf{R}^n : f_1(x_1, \ldots, x_n) = 0, \ldots, f_m(x_1, \ldots, x_n) = 0 \ \} \tag{22}$$

is called a variety. There is an inverse inclusion of variety and generating sets:

$$V(f_1, \ldots, f_m) \subset V(f_1, \ldots, f_{m-1}) \subset \cdots \subset V(f_1) \quad \text{and} \quad \{ f_1 \} \subset \{ f_1, f_2 \} \subset \cdots \subset \{ f_1, \ldots, f_m \}.$$

Also note that different subsets of polynomials can have the same variety. The zeroes Theorem addresses this very phenomenon.

***Theorem 29.*** (Hilbert 1888?)   Let $V(f)$ and $V(f_1, \ldots, f_m)$ be two varieties. If $V(f) = V(f_1, \ldots, f_m)$ then $f(x)^v = a_1 f_1(x) + \cdots + a_m f_m(x)$ for some integer $v \geq 1$.

The zeroes theorem gives a one-to-one correspondence between varieties and radical ideals in the ring of polynomials in $\mathbf{R}[x_1, \ldots, x_n]$.

***Theorem 30.***   Let $f_1, \ldots, f_n$ be relatively prime polynomials in $\mathbf{R}[x_1, \ldots, x_n]$. Then

(i) $V(f_1) \cap \cdots \cap V(f_n)$ is finite. In other words, the system of polynomial equations $f_1(x_1, \ldots, x_n) = 0, \ldots, f_n(x_1, \ldots, x_n) = 0$ has finitely many zeroes.

(ii) The $\mathbf{R}$- algebra $\mathbf{R}[x_1 \ldots x_n]/(f_1, \ldots, f_m)$ is finite dimensional.

Proof: The case $n = 2$ is covered in [K, p. 7].  ∎

***Theorem 31.***   (Bezout 1783)  If a system of polynomial equations $f_1(x_1, \ldots, x_n) = 0, \ldots, f_m(x_1, \ldots, x_n) = 0$ has finitely many zeroes, then it has at most $\deg(f_1) \times \deg(f_2) \times \cdots \times \deg(f_m)$ zeroes, counting multiplicities.

A system of polynomial equations has finitely many zeroes if and only if the corresponding polynomials do not have a nonconstant factor in common. Under this condition, the solution set

$$V(f_1, \ldots, f_m) = \{ (x_1, \ldots, x_n) \in \mathbf{R}^n : f_1(x_1, \ldots, x_n) = 0, \ldots, f_m(x_1, \ldots, x_n) = 0 \} \tag{23}$$

is a variety of zero dimension. The statement of this result holds exactly in the complex projective plane $P^{n+1}(\mathbb{C})$, but not in the affine real plane $P^n(\mathbb{R})$. In the case of a linear system of equations $f_1(x_1, \ldots, x_n) = 0, \ldots, f_m(x_1, \ldots, x_n) = 0$, of nonzero determinant, exempli gratia, $n = m$, the number of solutions is precisely $\deg(f_1) \times \deg(f_2) \times \cdots \times \deg(f_m) = 1$. Many specific results that give exact enumeration of the zeroes of zero dimension varieties are known, see the literature for details. A general result for the number of positive roots is the following.

***Theorem 32.***   (Khovanskii, 1999)   A system of $n$ polynomials equations in $n$ variables and with a total of $m$ monomials has at most $2^{m(m-1)/2}(n+1)^m$ positive roots.

***Theorem 33.***   A list of polynomials $f_1(x_1, \ldots, x_n), \ldots, f_m(x_1, \ldots, x_n)$ in the function field $\mathbf{R}[x_1 \ldots x_n]/(I(V))$ are algebraically independent over the ring $\mathbf{R}$ if there does not exist a nonconstant polynomial $g(x_1, \ldots, x_n)$ such that $g(f_1, \ldots, f_n) = 0$ in $\mathbf{R}[x_1, \ldots, x_n]/(I(V)$.

A functions field is the collection of rational functions defined by

$$\mathbf{R}[x_1, \ldots, x_n]/(I(V)) = \{ f(x)/g(x) : g(x) \neq 0, x \in V(I) \}. \tag{24}$$

***Theorem 34.***   The dimension of an affine variety $V(I)$ is the same as the maximal number of algebraically independent polynomial functions $f_i(x_1, \ldots, x_n)$ in the function field $\mathbf{R}[x_1, \ldots, x_n]/(I(V))$.

**Parametrizations Of A Varieties**

A parametrization of a variety is map $\phi : \mathbf{R}^m \rightarrow V \subset \mathbf{R}^n$ defined by the rational map $\phi(t_1, \ldots, t_m) = (\phi_1(t_1, \ldots, t_m), \ldots, \phi_n(t_1, \ldots, t_m))$.

As an example, the parametric form of the ellipse $x^2 + dy^2 = 1$ in $\mathbb{R}^2$ is described by the rational map $\phi(t) = \left( \dfrac{t^2 - d}{t^2 + d}, \dfrac{2t}{t^2 + d} \right)$, and the parametric form of the $n$th unit sphere $S^n = \{ (x_1, ..., x_n) : x_1^2 + \cdots + x_n^2 = 1 \}$ in $\mathbb{R}^n$ is described by the rational map

$$\phi(t_1, ..., t_{n-1}) = \left( \frac{2t_1}{t_1^2 + \cdots + t_{n-1}^2 + 1}, ..., \frac{2t_{n-1}}{t_1^2 + \cdots + t_{n-1}^2 + 1}, \frac{t_1^2 + \cdots + t_{n-1}^2 - 1}{t_1^2 + \cdots + t_{n-1}^2 + 1} \right). \tag{25}$$

### 3.2 Resultants

Let $f(x_1, ..., x_n) = a_k x_1^k + a_{k-1} x_1^{k-1} + \cdots + a_1 x_1 + a_0$, $g(x_1, ..., x_n) = b_m x_1^m + b_{m-1} x_1^{m-1} + \cdots + b_1 x_1 + b_0$, where $a_i$, $b_i \in \mathbf{R}[x_2, ..., x_n]$. The (Euler) *resultant matrix* $A = A(f, g)$ is defined by a $(k + m) \times (m + k)$ semi-circulant matrix. The matrix $A$ has $m$ columns of cyclic shifts of the coefficients of $f(x) = f(x_1, ..., x_n)$, and $k$ columns of cyclic shifts of the coefficients of $g(x) = g(x_1, ..., x_n)$.

The resultant matrix is derived from the linear system of equations

$$f(x) = 0, \; xf(x) = 0, \; x^2 f(x) = 0, \; ..., \; x^{m-1} f(x) = 0, \; g(x) = 0, \; xg(x) = 0, \; x^2 g(x) = 0, \; ..., \; x^{k-1} g(x) = 0. \tag{26}$$

There are $k + m$ linear equations in the variables $x_1 = x$, $x_2 = x^2$, ..., $x_k = x^k$, ..., $x_{k+m} = x^{k+m}$. A singular matrix $A(f, g)$ implies that $\gcd(f(x), g(x)) = $ polynomial of degree $> 1$. Otherwise, the matrix is nonsingular, and $\gcd(f(x), g(x)) = $ constant $\neq 0$.

**Definition 35.** The *resultant* of a pair of polynomials $f(x)$, $g(x) \in \mathbf{R}[x_1, ..., x_n]$ is defined as the determinant of the corresponding matrix, id est, $\mathrm{Res}(f, g, x_1) = \det(A) \in \mathbf{R}[x_2, ..., x_n]$.

A pair of relatively prime polynomials $f(x_1, ..., x_n)$, $g(x_1, ..., x_n) \in \mathbf{R}[x_1, ..., x_n]$ has $n$ (distinct) resultants, one for each variable $x_i$, $i = 1, 2, ..., n$. In particular,

(1) $r_1(x_2, x_3, ..., x_n) = \mathrm{Res}(g, f, x_1)$ effectively eliminates the variable $x_1$,
(2) $r_2(x_1, x_3, ..., x_n) = \mathrm{Res}(g, f, x_2)$ effectively eliminates the variable $x_2$,
...
(n) $r_n(x_1, x_2, ..., x_{n-1}) = \mathrm{Res}(g, f, x_n)$ effectively eliminates the variable $x_n$,

in the system of polynomial equations $f(x_1, ..., x_n) = 0$, $g(x_1, ..., x_n) = 0$.

Some useful and elementary properties of the resultants are stated below.

**Theorem 36.** Let $f(x_1, ..., x_n)$, $g(x_1, ..., x_n) \in \mathbf{R}[x_1, ..., x_n]$ be nonconstant polynomials. Then
(i) The resultant $\mathrm{Res}(f, g, x_1) \in \mathbf{R}[x_2, ..., x_n]$.
(ii) The reverse resultant $\mathrm{Res}(g, f, x_1) = (-1)^{km} \mathrm{Res}(f, g, x_1)$.
(iii) Multiplicative property $\mathrm{Res}(g, f, x_1) = \mathrm{Res}(f_1, g, x_1) \mathrm{Res}(f_2, g, x_1)$ if $f = f_1 f_2$.
(iv) $\mathrm{Res}(f, g, x_1) = f(x_1, ..., x_n)s + g(x_1, ..., x_n)t$, where $s$, $t \in \mathbf{R}[x_1, ..., x_n]$.

The formulas of the resultants are derived from the roots or the coefficients of the polynomials under consideration. Both techniques are very useful in theory and applications of the resultants.

*Coefficients Formulas.*
**Theorem 37.**  Let $f(x), g(x) \in \mathbf{R}[x_1,\ldots,x_n]$. The resultant is determinant of the corresponding matrix:
(i) $\mathrm{Res}(f, g, x_1) = \det(A) \in \mathbf{R}[x_2,\ldots,x_n]$,
(ii) $\mathrm{Re}\,s(f, g, x_1) = a_0^m \det(g(C_f) = (-1)^{km} b_0^k \det(f(C_g)$,

where $C_f$ and $C_g$ are the companion matrices of $f$ and $g$.

*Roots Formulas.*
**Theorem 38.**   Let $\alpha_1, \ldots, \alpha_k$, and $\beta_1, \ldots, \beta_m \in \mathbf{R}(x_2,\ldots,x_n)$ be the roots of the polynomials $f(x_1,\ldots,x_n)$ and $g(x_1,\ldots,x_n)$. Then

(i) $\mathrm{Re}\,s(f, g, x_1) = a_k^m b_m^k \prod_{i=1}^{k} \prod_{j=1}^{m} (\alpha_i - \beta_j) = a_k^m \prod_{i=1}^{k} g(\alpha_i) = (-1)^{km} b_m^k \prod_{i=1}^{m} f(\beta_i)$,

(ii) $\mathrm{Re}\,s(f, g, x_1) = b_m^{\deg(f) - \deg(g)} \prod_{i=1}^{t} r(\theta_i) = b_m^{\deg(f) - \deg(g)} \mathrm{Re}\,s(f, r, x_1)$,

where $f(x) = q(x)g(x) + r(x)$.

The roots of the polynomials $f$ and $g$ are assumed to be rational functions in the quotient fields or the algebraic closure. These expressions are of both theoretical and practical interests. Successive application of (ii) above reduces the size of the resultant matrix, thereby reducing the time complexity of computing $\mathrm{Res}(f(x), g(x), x)$. Another practical way of computing the resultant of a pair of relatively prime polynomials is by means of the Euclidean algorithm.

The root formula of the resultant is related to the determinant of the Vandermonde matrix $V = \left( \alpha_i^j \right), 1 \le i, j \le n$. Thus similar computation technique can be applied to both problems.

*Irreducibility Of The Resultants*
**Theorem 39.**   There exists a unique (up to sign) irreducible polynomial $\mathrm{Res}(f(x), g(x), x)$ in $\mathbf{R}[t_1,\ldots,t_d]$, $d = k + m$, which vanishes at the $(t_1,\ldots,t_d) = (a_1,\ldots,a_k, b_1,\ldots,b_m)$ whenever the polynomials $f(x) = a_k x^k + \cdots + a_1 x + a_0$ and $g(x) = b_m x^m + \cdots + b_1 x + b_0$ have a common root $x = (x_1,\ldots,x_n)$.

A proof of this result from the point of view of multivariable polynomials appears in [SL]. The irreducibility claim is valid for monic polynomials only, some other mild conditions might be required. This exception can be observed on the roots formula.

**Definition 40.**   A subset of polynomials $\{ f_1(x_1,\ldots,x_n), \ldots, f_m(x_1,\ldots,x_n) \} \subset \mathbf{R}[x_1,\ldots,x_n]$ are *algebraically independent* if and only if $\mathrm{Res}(f_i, f_j, x_k) \ne$ constant for all pairs $i \ne j$, and $k = 1, 2, \ldots, n$.

Relatively prime and algebraically independent polynomials are important in the elimination of all the variables but one variable in systems of polynomials equations.
**Discriminant**
The *discriminant* of a monic polynomial is the expression $\mathrm{disc}(f) = (-1)^{k(k-1)/2} \mathrm{Res}(f, f', x_1)$, where $f'$ is the derivative of $f$. The discriminant is widely used to identify polynomials with multiple roots. Similarly, the vanishing resultant $\mathrm{Res}(f, g, x_1) = 0$ implies that the polynomials $f$ and $g$ have a common factor of degree $\ge 1$.

*Multiplicative Property of the Discriminant*
(i) $\mathrm{disc}(fg) = \mathrm{disc}(f)\mathrm{disc}(g)\mathrm{Res}(f, g, x_1)^2$,
(ii) $\mathrm{disc}(fgh) = \mathrm{disc}(f)\mathrm{disc}(g)\mathrm{disc}(h)\mathrm{Res}(f, g, x_1)^2 \mathrm{Res}(f, h, x_1)^2 \mathrm{Res}(g\,h, x_1)^2$.

The multiplicative property follows from the definition and separation of the product into multiple products, one for each set of differences of roots, and the products of the mixed differences of roots.

The resultant method of elimination has several advantages over the related concept of Groebner bases method of elimination depending on the systems of equations:

(i) Order-free calculations: there is no need to have any order on the ring of polynomials.

(ii) Speed: there are several ways of computing the resultants effectively.

### 3.3 Groebner Bases

The extended Euclidean algorithm in the ring $\mathbf{R}[x]$ of polynomial in one variable gives a unique principal ideal (basis) $f(x) = a_1(x)f_1(x) + \cdots + a_m(x)f_m(x)$ of the ideal generated by the subset of polynomials $\{ f_1(x), \ldots, f_m(x) \} \subset \mathbf{R}[x]$. A Grobner basis is a generalization of this concept. The special basis of the ideal generated by the subset of polynomials $\{ f_1(x_1,\ldots,x_n), \ldots, f_m(x_1,\ldots,x_n) \} \subset \mathbf{R}[x_1,\ldots,x_n]$ is called a Grobner basis. The difficulty in the generalization arises from the fact that the ring $\mathbf{R}[x_1,\ldots,x_n]$ is not a principal ideal domain. This short note is intended to provide a rough working knowledge of Groebner bases, for finer details the reader should confer the literature, [CO], [BW], et cetera.

An *ordering* on the lattice $\mathbb{N}^n = \{ (\alpha_1, \ldots, \alpha_n) : \alpha_i \in \mathbb{N} \}$ of integers vectors $\alpha = (\alpha_1, \ldots, \alpha_n)$ is a relation on the set.

***Proposition* 41.** The lattice $\mathbb{N}^n$ has the following properties.

(i) $\mathbf{0} \preccurlyeq \alpha$, where $\mathbf{0} = (0, \ldots, 0)$ is the zero vector.

(ii) It is well ordered set: it has a least element with respect to an ordering.

(iii) It is an additive set: $\alpha \preccurlyeq \beta$ implies that $\alpha + \gamma \preccurlyeq \beta + \gamma$ for all integer vectors $\alpha, \beta, \gamma \in \mathbb{N}^n$.

The most common orderings are the followings:

1. Lexicographical ordering: $x^\alpha \preccurlyeq x^\beta$ if and only if $\alpha = \beta$ or there is an index such that the component $\alpha_j = \beta_j$ for $1 \leq j \leq i - 1$, and $\alpha_i < \beta_i$ for $i > j$.

This can be rephrased as $x^\alpha \preccurlyeq x^\beta$ if and only if $\alpha_i < \beta_i$ for the first index $i$ such that $\alpha_i \neq \beta_i$.

2. Inverse Lexicographical ordering: $x^\alpha \preccurlyeq x^\beta$ if and only if $\alpha = \beta$ or there is an index such that $\alpha_j = \beta_j$ for $i + 1 \leq j \leq n$ and $\alpha_i < \beta_i$ for $i > j$.

This can be rephrased as $x^\alpha \preccurlyeq x^\beta$ if and only if $\alpha_i < \beta_i$ for the last index $i$ such that $\alpha_i \neq \beta_i$.

3. Total degree lexicographical ordering: $x^\alpha \preccurlyeq x^\beta$ if and only if $|\alpha| = \alpha_1 + \cdots + \alpha_n < |\beta| = \beta_1 + \cdots + \beta_n$ or $|\alpha| = |\beta|$ and $\alpha_i < \beta_i$ for the first index $i$ such that $\alpha_i \neq \beta_i$.

4. Total degree inverse lexicographical ordering: $x^\alpha \preccurlyeq x^\beta$ if and only if $|\alpha| = \alpha_1 + \cdots + \alpha_n < |\beta| = \beta_1 + \cdots + \beta_n$ or $|\alpha| = |\beta|$ and $\alpha_i < \beta_i$ for the last index $i$ such that $\alpha_i \neq \beta_i$.

The ordering on the lattice $\mathbb{N}^n$ induces an ordering on the ring of polynomials $\mathbf{R}[x_1,\ldots,x_n]$. Under this scheme, there is a morphism $\alpha \rightarrow x^\alpha$. That is, each integer vector $\alpha = (\alpha_1, \ldots, \alpha_n) \in \mathbb{N}^n$ is mapped to a monomial $x^\alpha = x_1^{\alpha_1} x_2^{\alpha_2} \cdots x_n^{\alpha_n} \in \mathbf{R}[x_1,\ldots,x_n]$.

The time and space complexities of Groebner bases calculations greatly depend on the ordering on the ring of polynomials or equivalently the lattice $\mathbb{N}^n$. It is known that the reverse lexicographical ordering $1 \preccurlyeq x_1 \preccurlyeq x_2 \preccurlyeq \cdots \preccurlyeq x_n$ has lower complexity than lexicographical ordering $1 \succcurlyeq x_1 \succcurlyeq x_2 \succcurlyeq \cdots \succcurlyeq x_n$.

The unique leading term $LT(f) = x^\alpha$ of a polynomial $f(x) = a_\alpha x^\alpha + a_\beta x^\beta + \cdots + a_\mu x^\mu$ with respect to an ordering $\succcurlyeq$ is the monomial that satisfies the relation $x^\alpha \succcurlyeq x^\beta \succcurlyeq \cdots \succcurlyeq x^\gamma$.

A *Groebner* basis $G = \{\ g_1, \ldots, g_t\ \}$ of an ideal $I = <f_1, \ldots, f_t>$ is a basis such that $\{\ LT(g_1), \ldots, LT(g_t)\ \}$ generates the ideal .

The first elimination ideal $I_1$ eliminates the first variable $x_1$ in the Groebner basis, the second elimination ideal $I_2$ eliminates the first variable $x_2$ in the Groebner basis, and so on.

***Theorem 42.*** (?)   Let $I = <f_1, \ldots, f_m>$ be an ideal in $\mathbf{R}[x_1, \ldots, x_n]$ , and let $G = \{\ g_1, \ldots, g_k\ \}$ be a Groebner basis for $I$. Then every polynomial $f(x) \in \mathbf{R}[x_1, \ldots, x_n]$ has a representation as $f(x) = a_1 f_1(x) + \cdots + a_m f_m(x) + r(x)$, where $a_i(x), r(x) \in \mathbf{R}[x_1, \ldots, x_n]$, and the remainder $r(x)$ is unique.

Given an ideal $I = <f_1, \ldots, f_m>$, the $k$th elimination ideal $I_k = I \cap \mathrm{R}[x_{k+1}, \ldots, x_n]$ is the smaller ideal generated by the subset of $n - k$ polynomials $\{\ h_1, \ldots, h_{n-k}\ \} \subset \{ f_1, \ldots, f_m \ \}$.

***Theorem 43.*** (Elimination theorem)   Let $G = \{\ g_1, \ldots, g_t\ \}$ be a Groebner basis of the ideal $I = <f_1, \ldots, f_m>$ in $\mathbf{R}[x_1, \ldots, x_n]$ with respect to the ordering $x_1 \succcurlyeq x_2 \succcurlyeq \cdots \succcurlyeq x_n$. Then the set $G_k = G \cap \mathbf{R}[x_1, \ldots, x_n]$ is a Groebner basis of the $k$th elimination ideal $I_k$.

Proof: See [CO, p. 115], [BW, p.399].

***Theorem 44.*** A system of polynomial equations $f_1 = 0, \ldots, f_m = 0$ is not solvable
(i) If the corresponding Groebner basis $G = \{\ 1\ \}$.
(ii) If $a_1 f_1(x) + \cdots + a_m f_m(x) = 1$ for some $a_i(x) \in \mathbf{R}[x_1, \ldots, x_n]$.

### 3.4 Construction Of Algebraically Independent Polynomials

The application of lattice reduction theory to the theory of polynomial equations is a recent phenomenon. The application to linear integer equations linear modular equations has an earlier beginning, see [BJ] and [LR], and application to nonlinear algebraic equations seems to have been the works of [HD] and [VE]. Later the technique for nonlinear congruent equations was significantly improved and extended to integer equations in [CR]. Specifically, the range of the roots $x$ of the polynomial equation $f(x) = a_d x^d + \cdots + a_1 x + a_0$ modulo $N$ that can be determined in deterministic logarithmic time complexity was extended from $|x| \leq N^{2/(d(d+1))}$ to $|x| \leq N^{1/d-\varepsilon}$, $\varepsilon > 0$. This is accomplished by replacing the original basis of the polynomials lattice.

Since then over a dozen papers and a few dissertations (see [D], [M], and [J]), have been published on the applications of lattice reductions theory to integer factorization and cryptography. These more recent works have simplified the theory and its practical aspect. The recent advances (as the improved larger ranges of the roots of polynomial equations) in polynomials equations and lattice reduction methods are mostly direct results of the new streamlined polynomial bases of the polynomial lattices, see [ER], [BA], [BM], [J].

The $l_p$ norm of a polynomial is defined by the usual expression

$$\| f(x) \|_p = \sqrt[p]{\sum | a_\alpha |^p} \ , \tag{27}$$

where $0 < p \le \infty$. The cases $p = 2$ and $p = \infty$ are widely used in the analysis of algorithms of polynomials. These are supremum norm and standard norm given by the relations

$$\| f(x) \|_\infty = \max \{ | a_\alpha | : \alpha \ne 0 \} \ \text{ and } \ \| f(x) \|_2 = \sqrt{\sum | a_\alpha |^2} \tag{28}$$

respectively. The Lehmer measure of a multivariable polynomial is defined by

$$M(f) = \exp \int_0^1 \cdots \int_0^1 \ln \left| f(e^{i2\pi s_1}, ..., e^{i2\pi s_n}) \right| ds_1 \cdots ds_n . \tag{29}$$

***Theorem* 45.** (i) Let $f(x_1, ..., x_n)$ be of maximum degree $n_i$ in the variable $x_i$, and let $\alpha = (\alpha_1, ..., \alpha_n)$ be a vector of positive integers. Then

(i) $|a_\alpha| \le \binom{n_1}{\alpha_1} \binom{n_2}{\alpha_2} \cdots \binom{n_n}{\alpha_n} M(f)$ and (ii) $H(f) \le 2^{n_1 + n_2 + \cdots n_n - n} M(f)$.

A proof of this result is given in [PV, p. 149].

***Theorem* 46.** ([ST]) Let $f(x_1, ..., x_n)$ and $g(x_1, ..., x_n)$ be two non-zero polynomials over $\mathbb{Z}$ of maximum degree $d$ in each variable separately such that $g(x_1, ..., x_n)$ is a multiple of $f(x_1, ..., x_n)$ in $\mathbb{Z}[x_1, ..., x_n]$. Then
$$\| g \|_2 \ge 2^{-(d+1)^n + 1} \| f \|_\infty$$

***Theorem* 47.** ([NH]) Let $f(x_1, ..., x_n) \in \mathbb{Z}[x_1, ..., x_n]$ be a polynomial in $n$-variables with $w > 0$ nonzero terms. Suppose that $f(x_1, ..., x_n) \equiv 0 \bmod N$, and $\| f(x_1 X_1, ..., x_n X_n) \|_\infty < N w^{-1/2}$, where $| x_i | < X_i$ for $i = 1, 2, ..., n$. Then $f(x_1, ..., x_n) = 0$ holds over the integers $\mathbb{Z}$.

***Theorem* 48.** ([CR]) Let $f(x, y) \in \mathbb{Z}[x, y]$ be irreducible with maximum degree $d$ in $x, y$ separately. Let $| X |$, $| Y |$ be upper bounds on the desired integer solution $(x_0, y_0)$ and let $W = \max \{ | a_{i,j} XY | \mid 0 \le i, j \le d \}$. If $XY \le W^{2/3d}$ then all integer pairs $(x_0, y_0)$ such that $f(x_0, y_0) = 0$, $| x_0 | \le X$ and $| y_0 | \le Y$ can be found in time polynomial in log $W$ and $2^d$.

Although Theorem 47 calls for an irreducible polynomial, Theorem 45 seems to imply that reducible polynomials work as well as irreducible polynomials whenever $\| g \|_2 < 2^{-(d+1)^n + 1} \| f \|_\infty$ holds.

***Definition* 49.** Let $f(x_1, ..., x_n) \in \mathbb{Z}[x_1, ..., x_n]$ and let $S, T \subset \mathbb{Z}[x_1, ..., x_n]$ be subsets of monomials $x^\alpha = x_1^{\alpha_1} x_2^{\alpha_2} \cdots x_n^{\alpha_n}$, where $\alpha = (\alpha_1, ..., \alpha_n)$ is an integer vector with $\alpha_i \ge 0$. The monomials sets $S$ and $T$ are called admissible for $f(x) = f(x_1, ..., x_n)$ if and only if
(i) For every monomial $x^\alpha \in S$, the polynomial $x^\alpha f(x)$ is defined over $T$.
(ii) For every polynomial $g(x) = g(x_1, ..., x_n)$ defined over $T$, if $g(x) = f(x)h(x)$ for some $h(x)$, then $h(x)$ is defined over $S$.

The 2-variable, 3-variable, and 4-variable versions of this definition are discussed in detailed in [BM], [BA], and [JM] respectively.

The set $S$ consists of the shift monomials $\{\, x^\beta = x_1^{\beta_1} x_2^{\beta_2} \cdots x_n^{\beta_n} : \beta \in B \,\}$, where $B$ is some index set. The shift monomials sets $S$ is used to generate a partial basis $\{\, x^\beta f(x) : \beta \in B \,\}$ of the polynomial lattice. The selection of an optimum set $S$, which yields the widest ranges of the variables $|\,x_1\,| < X_1, \ldots, |\,x_n\,| < X_n$ is an open problem, see [BM].

The ranges $|\,x_1\,| < X_1, \ldots, |\,x_n\,| < X_n$ of the variables are maximal on small Newton polygons. Exempli gratia, the early result in [CR] for polynomials $f(x, y)$ in 2-variable has the widest ranges of the variables $|\,x\,| < X$ and $|\,y\,| < Y$ if the Newton polygon of the polynomial $f(x, y)$ is a triangle than for a rectangle.

A Newton polygon is the convex set enclosed by the convex hull, and the convex hull is the set of indices $\mathrm{conv}(f) = \{\, \alpha : a_\alpha \neq 0 \text{ is a coefficient of } f(x) \,\}$.

The set of monomials $T$ consists of the monomials in the expansion of $x^\beta f(x_1, \ldots, x_n)$ over some set of indices $\beta \in B$. In many cases $T = S + \mathrm{ind}(f)$.

For example, in the three variables case, a polynomial $f(x, y, z)$ is defined over a subset of monomials $S = \{\, x^\beta = x^{\beta_1} y^{\beta_2} z^{\beta_3} : \beta = (\beta_1, \beta_2, \beta_3) \in B \,\}$ if $f(x, y, z)$ can be written as linear combination of monomials in $S$. Similarly, the polynomial $h(x, y, z)$ is defined over $T$, if $h(x, y, z) = g(x, y, z) f(x, y, z)$ for some polynomial $g(x, y, z)$. The ordered pair $(S, T)$ is said to be admissible for $f$, see [BM], [BA] for detailed discussions.

The integer $d_1$, $d_2$ and $d_3$ are the maximum degree of the polynomial $f(x, y, z)$ in the variables $x$, $y$, and $z$ respectively. Furthermore, the integers $s_1$, $s_2$ and $s_3$ are defined by the sums

$$ s_1 = \sum_{(i,j,k) \in T \setminus S} i \quad , \quad s_2 = \sum_{(i,j,k) \in T \setminus S} j \quad , \quad \text{and} \quad s_3 = \sum_{(i,j,k) \in T \setminus S} k \quad . \tag{30}$$

Similarly, the height and norm of a polynomial $f(x, y, z) = \sum_{0 \le i,j,l \le d} a_{i,j,l}\, x^i y^j z^l \in \mathbb{Z}[x,y,z]$ of maximum degree $\deg(f) = d$ in the variables $x$ and $y$ are given by the expressions $\|f(x,y,z)\|_\infty = \max\{\, |a_{0,0,0}|, |a_{0,0,1}|, \ldots, |a_{d,d,d}|\, \}$ and $\|f(x, y, z)\|_2 = \sqrt{a_{0,0,0}^2 + a_{0,0,1}^2 + \cdots + a_{d,d,d}^2}$ respectively.

Let $f(x, y, z)$ be an irreducible polynomial in $\mathbb{Z}[x, y, z]$ of height $\|f(xX, yY, zZ)\|_\infty = W$, and let $(x_0, y_0, z_0)$ be a small integer root such that $|\,x_0\,| < X$, $|\,y_0\,| < Y$ and $|\,z_0\,| < Z$.

**Theorem 50.** ([BA]) If $S$ and $T$ are admissible sets for $f(x, y, z)$, then an algebraically independent polynomial $g(x, y, z)$ which has $(x_0, y_0, z_0)$ as a root over the integers can be found in logarithm time, provided that

$$ X^{s_1} Y^{s_2} Z^{s_3} < W^s\, 2^{-(6+c)(d_1^2 + d_2^2 + d_3^2)s} , \tag{31}$$

where it is assumed that $(t - s)^2 \le cs(d_1^2 + d_2^2 + d_3^2)$ for some constant $c$.

This is a general result for arbitrary polynomials in three variables. A few specialized cases of polynomials in three and four variables have been worked out in details. Two such cases are considered here.

**Theorem 51.** ([ER]) Let $f(x, y, z) = c_4 xy + c_3 x + c_2 y + c_1 z + c_0 \in \mathbb{Z}[x, y, z]$ be an irreducible polynomial of height $\|f(xX, yY, zZ)\|_\infty = W$ and with a small integer root $(x_0, y_0, z_0)$ such that $|\,x_0\,| < X$, $|\,y_0\,| < Y$ and $|\,z_0\,| < Z$. Suppose that the inequality

$$ X^{3+3\tau} Y^{3+6\tau+3\tau^2} Z^{2+3\tau} < W^{2+3\tau-\varepsilon} , \tag{32}$$

where $\tau > 0$ is a lattice parameter, holds. Then there exists a pair of linearly independent polynomials $f_1(x, y, z)$ and $f_2(x, y, z)$ not multiple of $f(x, y, z)$, with a common root. Furthermore, the polynomials are generated in deterministic logarithm time complexity.

The complete analysis of this and other special cases and the corresponding polynomials bases of the polynomials lattices are given in [ER], [JM].

The two polynomials $f_1(x, y, z)$ and $f_2(x, y, z)$ are not guaranteed to be algebraically independent. However, the two pairs of polynomials $f(x, y, z)$, $f_1(x, y, z)$ and $f(x, y, z)$, $f_2(x, y, z)$ are guaranteed to be algebraically independent. Recent advances in the construction of three algebraically independent polynomials are discussed in [BA].

***Corollary* 52.**    Let $f(x, y, z) \in \mathbb{Z}[x,y,z]$ be an irreducible polynomial of maximal degree $d > 0$ in $x$, $y$, and $z$, (or of total degree $d > 0$ in $x$, $y$, and $z$) and let $\| f(xX, yY, zZ) \| = W$ be the height of $f$. Suppose that there exists a triple $(x_0, y_0, z_0)$ such that $f(x_0, y_0, z_0) = 0$, where $0 \leq | x_0 y_0 z_0 | < W^{2/3d}$, and $0 \leq | z | \leq O((\log N)^B)$, $B > 0$ constant. Then the solution $(x_0, y_0, z_0)$ can be determined in deterministic logarithmic time.